\documentclass[11pt]{article}
\usepackage{amsmath}
\usepackage{pxfonts}
\usepackage{yfonts}

\usepackage{dsfont}
\usepackage{marvosym}
\usepackage{graphicx}
\usepackage{relsize}

\def\bel{\begin{equation}\label}
\def\eeq{\end{equation}}
\def\ds{\displaystyle}
\def\endproof{\hphantom{MM}
\hfill\llap{$\square$}\goodbreak}

\def\mt{\longrightarrow}
\def\v{\vskip 1em}
\def\vsk{\vskip 40em}
\def\ve{\varepsilon}
\def\R{\mathbb R}
\def\Z{\mathbb Z}
\def\C{\mathfrak{C}}

\def\N{{\bf N}}

\def\L{{\bf L}}

\def\bar{\overline}

\def\I{{\bf I}}

\def\M{{\bf M}}

\def\Cup{{\bigcup}}
\def\Cap{{\bigcap}}

\def\alpha{\alphaup}
\def\beta{\betaup}
\def\gamma{\gammaup}
\def\delta{\deltaup}
\def\xi{{\xiup}}
\def\eta{{\etaup}}
\def\tau{{\tauup}}
\def\rho{{\rhoup}}
\def\phi{{\phiup}}
\def\psi{{\psiup}}
\def\lambda{{\lambdaup}}
\def\omega{\omegaup}
\def\varphi{{\varphiup}}
\def\gamma{{\gammaup}}

\newtheorem{prop}{Proposition}[section]
\newtheorem{remark}{Remark}[section]

\begin{document}
 \[\begin{array}{cc}\hbox{\LARGE{\bf Fractional integration with singularity on Light-cone}}
 \\\\
 \hbox{\LARGE{\bf I: {\it natural setting}}} 
\end{array}\]

 \[\hbox{Zipeng Wang}\]
  \[\begin{array}{ccc}
 \hbox{ Department of Mathematics, Westlake university}\\
 \hbox{Cloud town, Hangzhou of China}
 \end{array}\]
 \begin{abstract}
 In this first part of our project, we prove a classical Hardy-Littlewood-Sobolev result for a new family of 
fractional integral operators whose kernel has singularity appeared on the light-cone in $\R^{n+1}$. \end{abstract}
  \section{Introduction}
 \setcounter{equation}{0}
In 1928, Hardy and Littlewood \cite{Hardy-Littlewood}  first established a regularity theorem on $\R$ for certain convolution operators. Ten years later, Sobolev \cite{Sobolev} extended this result to higher dimensions. Today, it is known as Hardy-Littlewood-Sobolev inequality for fractional integral operators. 
Their regularity  has been extensively studied over the past several decades, for example  by Stein and Weiss \cite{Stein-Weiss}, Ricci and Stein \cite{Ricci-Stein}, Strichartz \cite{Strichartz},  Fefferman and Muckenhoupt \cite{Fefferman-Muckenhoupt},   Muckenhoupt and Wheeden \cite{Muckenhoupt-Wheeden},  Sawyer and Wheeden \cite{Sawyer-Wheeden},   Perez \cite{Perez}, Sawyer and Wang \cite{Sawyer-Wang} and Wang \cite{Wang}. 

We consider
\bel{I_alpha}
\begin{array}{cc}\ds
 \Big(\I_{\alpha}f\Big)(x,y)~=~\iint_{\Lambda} f(x-y,t-s)\left({1\over |s|+|y|}\right)^{n-\alpha}\left({1\over |s|-|y|}\right)^{1-{\alpha\over n}}dy ds,
 \\\\ \ds
 \Lambda~=~\Big\{ (y,s)\in\R^n\times\R~\colon |s| >|y|\Big\}.
 \end{array} 
\eeq
Write $\C$ as a  generic constant.  Our main result is given below.

{\bf Theorem A~~}   {\it Let $\I_\alpha$ be defined in (\ref{I_alpha}) for $0<\alpha<n$.  We have
\bel{Result}
\left\| \I_{\alpha} f\right\|_{\L^q(\R^{n+1})}~\leq~\C_{p~q}~\left\| f\right\|_{\L^p(\R^{n+1})}
\eeq
for $1<p<q<\infty$ if and only if
\bel{Formula}
{\alpha\over n}~=~{1\over p}-{1\over q}.
\eeq
}

The necessity of (\ref{Formula}) can be verified by changing dilations, as a standard exercise. 

The paper is organized as follows. In  section 2, we prove the special case when $n=1$ by carrying out an iteration argument. This answers a question proposed by Oberlin \cite{Oberlin} for $n=1$.

In order to prove {\bf Theorem A} for $n\ge2$, we introduce a new framework  in section 3
where $\Lambda$ given in (\ref{I_alpha}) is decomposed into infinitely many discrete variant of dyadic cones, denoted by  $\Lambda_\ell,\ell\ge0$. 
Each one of them is a collection of truncated regions having the same eccentricity  $2^{-\ell}$ where we develop a pointwise estimate in the same sprit of Hedberg \cite{Hedberg}. 

In section 4, we  show that every partial operator defined on $\Lambda_\ell$ satisfies the desired regularity.  

Section 5 and 6 are devoted to some preliminary estimates. We prove a crucial lemma of almost orthogonality in section 7.

We assume $f\ge0$ in dealing with convolution operators of positive kernels.

\section{Iteration estimate on $\R^2$}
\setcounter{equation}{0}
Let $n=1$ in (\ref{I_alpha}). We have
 \bel{I_alpha R^2}
 \begin{array}{lr}\ds
 \Big(\I_{\alpha}f\Big)(x,y)~=~\iint_{\Lambda} f(x-y,t-s)\left({1\over |s|+|y|}\right)^{1-\alpha}\left({1\over |s|-|y|}\right)^{1-\alpha}dy ds 
 \\\\ \ds~~~~~~~~~~~~~~~~
 ~=~\iint_{\Lambda} f(x-y,t-s)\left({1\over s^2-|y|^2}\right)^{1-\alpha}dy ds  
 \\\\ \ds~~~~~~~~~~~~~~~~
 ~<~\iint_{\R^2} f(x-y,t-s)\left({1\over |y+s|}\right)^{1-\alpha}\left({1\over |y-s|}\right)^{1-\alpha}dy ds
 \end{array}
 \eeq
whose kernel is symmetric $w.r.t$ $(y+s,y-s)\in\R\times\R$.

By changing variables 
\bel{change variables}
z~=~{x+t\over 2},~~ w~=~{x-t\over 2},\qquad u~=~{y+s\over 2}, ~~v~=~{y-s\over 2}
\eeq
in (\ref{I_alpha R^2}), we have
\bel{I_alpha'}
\begin{array}{lr}\ds
~~~~~~~\iint_{\R^2} f(x-y,t-s)\left({1\over |y+s|}\right)^{1-\alpha}\left({1\over |y-s|}\right)^{1-\alpha}dy ds
\\\\ \ds
 ~=~2^{2\alpha-1}\iint_{\R^2} f(z-u+w-v,z-u-w+v)\left({1\over |u|}\right)^{1-\alpha}\left({1\over |v|}\right)^{1-\alpha}dudv.
 \end{array}
 \eeq
 Denote 
 \bel{f_z, F_w}
 \begin{array}{cc}\ds
 f_z(w-v)~=~f(z+w-v, z-w+v),
 \\\\ \ds
  F_w(z-u)~=~\int_{\R} f(z-u+w-v,z-u-w+v)\left({1\over |v|}\right)^{1-\alpha} dv
  \end{array}
 \eeq
 for $z,w\in\R$.
 
We recall a classical result obtained by Hardy, Littlewood and Sobolev  \cite{Hardy-Littlewood}-\cite{Sobolev}:

 \v
 {\bf Hardy-Littlewood-Sobolev theorem} ~~~{\it Let $0<\alpha<\N$. We have
 \bel{Fractional}
 \begin{array}{lr}\ds
 \left\{\int_{\R^\N} \left\{\int_{\R^\N} f(y)\left({1\over |x-y|}\right)^{\N-\alpha}dy\right\}^qdx\right\}^{1\over q}
 ~\leq~\C_{p~q}\left\{\int_{\R^\N} \Big(f(x)\Big)^pdx\right\}^{1\over p}
 \end{array}
 \eeq
 for $1<p<q<\infty$ if and only if
 \bel{FORMULA}
 {\alpha\over\N}~=~{1\over p}-{1\over q}.
 \eeq
 }
 
 Let $\alpha=1/p-1/q$.
 We have 
 \bel{Iteration}
 \begin{array}{lr}\ds
 \left\|\I_\alpha f\right\|_{\L^q(\R^2)}
 \\\\ \ds
 ~=~ \left\{\iint_{\R^2}  \left\{\iint_{\R^2} f(x-y,t-s)\left({1\over |y+s|}\right)^{1-\alpha}\left({1\over |y-s|}\right)^{1-\alpha}dy ds\right\}^q dxdy\right\}^{1\over q}
 \\\\ \ds
 ~=~ 2^{2\alpha}\left\{\iint_{\R^2}  \left\{\iint_{\R^2} f(z-u+w-v,z-u-w+v)\left({1\over |u|}\right)^{1-\alpha}\left({1\over |v|}\right)^{1-\alpha}dudv\right\}^q dzdw\right\}^{1\over q} 
 \\\\ \ds
 ~=~2^{2\alpha}\left\{\int_{\R}\left\{\int_{\R}\left\{\int_{\R} F_w(z-u)\left({1\over |u|}\right)^{1-\alpha}du \right\}^q dz\right\}dw\right\}^{1\over q}
 \\\\ \ds
 ~\leq~\C_{p~q} ~\left\{\int_{\R}\left\{\int_{\R} \Big(F_w(z)\Big)^p dz\right\}^{q\over p}dw\right\}^{1\over q} \qquad\hbox{\small{by Hardy-Littlewood-Sobolev theorem}}
 \\\\ \ds
 ~\leq~\C_{p~q} ~\left\{\int_{\R}\left\{\int_{\R} \Big(F_w(z)\Big)^q dw\right\}^{p\over q}dz\right\}^{1\over p} 
 \qquad \hbox{\small{by Minkowski integral inequality}}
 \\\\ \ds
 ~=~\C_{p~q} ~\left\{\int_{\R}\left\{\int_{\R} \left\{\int_{\R} f(z+w-v,z-w+v)\left({1\over |v|}\right)^{1-\alpha} dv\right\}^q dw\right\}^{p\over q}dz\right\}^{1\over p}  
 \\\\ \ds
 ~=~\C_{p~q} ~\left\{\int_{\R}\left\{\int_{\R} \left\{\int_{\R} f_z(w-v)\left({1\over |v|}\right)^{1-\alpha} dv\right\}^q dw\right\}^{p\over q}dz\right\}^{1\over p}   
 \\\\ \ds
 ~\leq~\C_{p~q}~\left\{\iint_{\R^2} \Big(f_z(w)\Big)^pdwdz\right\}^{1\over p}
 \qquad\hbox{\small{by Hardy-Littlewood-Sobolev theorem}} 
  \\\\ \ds
 ~=~\C_{p~q}~\left\{\iint_{\R^2} \Big(f(z+w,z-w)\Big)^pdzdw\right\}^{1\over p}
 \\\\ \ds
 ~=~\C_{p~q}~\left\| f\right\|_{\L^p(\R^2)}
 \qquad f\in\L^p(\R^2). 
 \end{array}
 \eeq

 \section{Cone decomposition on $\R^n\times\R$}
 \setcounter{equation}{0}
Let $\ell$ be an nonnegative integer. We define
 \bel{Partial}
 \Big(\Delta_\ell \I_{\alpha} f\Big)(x,t)~=~\iint_{\Lambda_\ell} f(x-y,t-s)\left({1\over |s|+|y|}\right)^{n-\alpha}\left({1\over |s|-|y| }\right)^{1-{\alpha\over n}}dy ds
\eeq
where
\bel{Cone}
\begin{array}{cc}\ds
\Lambda_\ell~=~\Cup_{j\in\mathbb{Z}}~ \Lambda_{\ell j},
\\\\ \ds
\Lambda_{\ell j}~=~\left\{(y,s)\in\R^n\times(0,\infty)~\colon 2^j\leq |s|+|y| < 2^{j+1}, 2^{j-\ell}\leq |s|-|y|<2^{j-\ell+1}\right\}.
\end{array}
\eeq 
We momentarily consider a partial sum operator $\I_\alpha^\eta~\doteq~\sum_{\ell=0}^\eta \Delta_\ell\I_\alpha $ for some $\eta\ge1$ of which $\I_\alpha=\lim_{\eta\mt\infty}\I_\alpha^\eta$. Our resulting estimate will be independent from its value. 
\begin{figure}[h]
\centering
\includegraphics[scale=0.40]{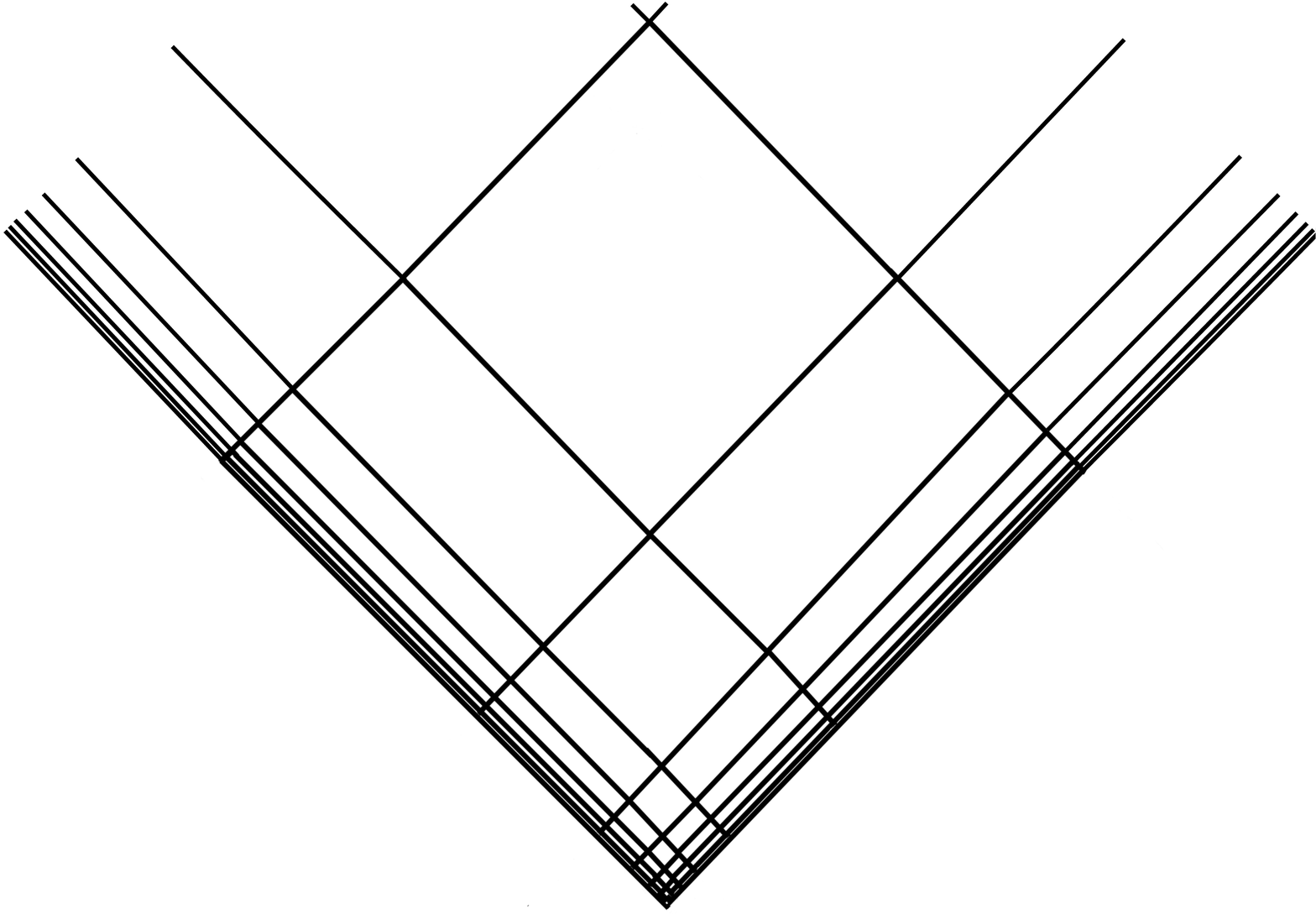}
\caption{ $\ds\Cup_{\ell\ge0, j\in\Z}\Lambda_{\ell j}$ inside  the upper half space where $s>0$ is in the vertical direction.}
\end{figure}

Let $n\ge2$ from now on.   
We consider a collection of points denoted by $\left\{y^\nu_\eta\right\}_\nu$  that are equally distributed on  the unit sphere $\mathds{S}^{n-1}\subset\R^n$ with grid length equal to $\mathfrak{B} 2^{-\eta}$ for some ${1\over 2}\leq\mathfrak{B}\leq 2$.  Note that there are at most $\C2^{\eta\left(n-1\right)}$ such elements . 
 
Define a geometric cone
 \bel{Gamma}
\Gamma^{\nu}_\eta~=~\left\{y\in\R^n~\colon~\left| {y\over |y|}-y^{\nu}_\eta\right|~\leq~2\times2^{1-\eta}\right\}
\eeq
whose central direction is $y^{\nu}_\eta$.
 
Observe that 
\bel{Inclusion R^v_j} 
\Lambda_{\ell j}\cap \Gamma^\nu_\eta~\subset~R^\nu_{\eta~\ell j}
\eeq
for some rectangle $R^\nu_{\eta~\ell j}\subset\R^{n}\times\R$ having  side lengths equal to a constant multiple of $ 2^j, 2^{j-\ell}$ and $ 2^{j-\eta}$ such that
\bel{R size}
\left| R^\nu_{\eta~\ell j}\right|~=~\C ~2^j 2^{j-\ell}2^{(j-\eta)(n-1)}
\eeq
for every $0\leq \ell\leq \eta$.

A maximal operator $\M^\nu_\eta$  defined by
\bel{M^v} 
 \Big(\M^\nu_\eta f\Big)(x,t)~=~\sup_{0\leq\ell\leq \eta,~ j\in\Z} ~{1\over 2^j2^{j-\ell}2^{(j-\eta)(n-1)}}\iint_{\Lambda_{\ell j}\cap\Gamma^\nu_\eta} f(x-y,t-s)dyds.
 \eeq
is bounded on $\L^p(\R^{n+1})$ for $1<p<\infty$. Recall that there are at most $\C 2^{\eta(n-1)}$ elements in $\{y^\nu_\eta\}_\nu$. 

Define
\bel{M*}
\Big(\M_\eta f\Big)(x,t) ~=~2^{-\eta(n-1)}\sum_\nu \Big(\M^\nu_\eta f\Big)(x,t).
\eeq
By Minkowski inequality, we clearly have $\M_\eta\colon\L^p(\R^{n+1})\mt\L^p(\R^{n+1}), 1<p<\infty$.

Let 
\bel{theta}
\vartheta_\ell(x,t)~=~{1\over\left\| f\right\|_{\L^p(\R^{n+1})}^p}\iint_{\Lambda_\ell}f^p(x-y,t-s)dyds
\eeq
for every $\ell\ge0$. It is clear that
\bel{theta sum}
\sum_{\ell=0}^\infty \vartheta_{\ell}(x,t)~\leq~1.
\eeq
We will show that for $a.e$ $(x,t)\in\R^n\times\R$,
  \bel{Regularity est}
 \begin{array}{lr}\ds
 \Big(\Delta_\ell \I_{\alpha} f\Big)(x,t)
 ~\leq~\C_{p~q}~ \vartheta_\ell^{{1\over p}-{1\over q}}(x,t) \Big(\M_\eta f\Big)^{p\over q}(x,t) \left\| f\right\|_{\L^p(\R^{n+1})}^{1-{p\over q}}
 \end{array}
\eeq
for every $0\leq \ell\leq \eta$.

Suppose $2p\leq q$. We give an heuristic estimate as follows. 
\bel{Heuristic}
 \begin{array}{lr}\ds
 \iint_{\R^{n+1}}  \sum_{\ell=0}^\eta \Big(\Delta_\ell \I_{\alpha} f\Big)^q(x,t) dxdt 
 \\\\ \ds
 ~\leq~\C_{p~q}~\left\| f\right\|_{\L^p(\R^{n+1})}^{q-p}\iint_{\R^{n+1}}\sum_{\ell=0}^\eta \vartheta_\ell^{{q\over p}-1}(x,t) \Big(\M_\eta f\Big)^{p}(x,t)  dxdt \qquad \hbox{\small{by (\ref{Regularity est})}}
  \\\\ \ds
 ~\leq~\C_{p~q}~\left\| f\right\|_{\L^p(\R^{n+1})}^{q-p}\iint_{\R^{n+1}}\left\{\sum_{\ell=0}^\infty \vartheta_\ell(x,t)\right\} \Big(\M_\eta f\Big)^{p}(x,t)  dxdt  
\\\\ \ds
  ~\leq~\C_{p~q}~\left\| f\right\|_{\L^p(\R^{n+1})}^{q-p}\iint_{\R^{n+1}} \Big(\M_\eta f\Big)^{p}(x,t)  dxdt \qquad\hbox{\small{by (\ref{theta})-(\ref{theta sum})}}
  \\\\ \ds
  ~\leq~\C_{p~q}~\left\| f\right\|_{\L^p(\R^{n+1})}^{q}. 
  \end{array}
 \eeq
 Note that the estimate in (\ref{Heuristic}) holds for every $\eta\ge1$. 
 
Let $h_i\in\Z, i=1,2,\ldots,q-1$ and $2\leq q\in\Z$. We have
\bel{Sum I_alpha}
\begin{array}{lr}\ds
\iint_{\R^{n+1}} \Big(\I_\alpha f\Big)^q(x,t)dxdt~=~
\\\\ \ds
\sum_{h_i,i=1,2,\ldots,q-1}\iint_{\R^{n+1}}\sum_{\ell\ge h_i, i=1,2,\ldots,q-1} \Big(\Delta_\ell \I_\alpha f\Big)(x,t) \prod_{i=1}^{q-1}\Big(\Delta_{\ell-h_i}\I_\alpha f\Big)(x,t) dxdt.
\end{array}
\eeq
After reordering the product $(\Delta_\ell \I_\alpha f) \prod_{i=1}^{q-1}(\Delta_{\ell-h_i}\I_\alpha f)$, it is suffice to consider $\ell\ge\ell-h_i$ for every $i=1,2,\ldots,q-1$ in the summand of (\ref{Sum I_alpha}).  
We therefore further require $h_i, i=1,2,\ldots,q-1$ to be nonnegative integers.

By applying H\"{o}lder inequality twice, we have
\bel{Ortho Expansion}
\begin{array}{lr}\ds
\iint_{\R^{n+1}}\sum_{\ell\ge h_i,i=1,2,\ldots,q-1} \Big(\Delta_\ell \I_\alpha f\Big)(x,t) \prod_{i=1}^{q-1}\Big(\Delta_{\ell-h_i}\I_\alpha f\Big)(x,t) dxdt
\\\\ \ds
~\leq~ \iint_{\R^{n+1}} \prod_{i=1}^{q-1} \left\{\sum_{\ell\ge h_i, i=1,2,\ldots,q-1} \Big(\Delta_\ell \I_\alpha f\Big)(x,t) \Big(\Delta_{\ell-h_i}\I_\alpha f\Big)^{q-1}(x,t)\right\}^{1\over q-1}dxdt
\\\\ \ds
~\leq~\prod_{i=1}^{q-1} \left\{\iint_{\R^{n+1}} \sum_{\ell\ge h_i, i=1,2,\ldots,q-1} \Big(\Delta_\ell \I_\alpha f\Big)(x,t) \Big(\Delta_{\ell-h_i}\I_\alpha f\Big)^{q-1}(x,t)dxdt\right\}^{1\over q-1}.
\end{array}
\eeq
We prove {\bf Theorem A} by showing that $\Delta_\ell\I_\alpha, \ell\ge0$ enjoy certain almost orthogonality. 
 \v

{\bf Lemma of almost orthogonality:} {\it Let ${\alpha\over n}={1\over p}-{1\over q}, 1<p<q<\infty$. There exists an $\ve=\ve(p,q)>0$ such that
\bel{Ortho Result}
\begin{array}{lr}\ds
\iint_{\R^{n+1}}\sum_{\ell= h}^\infty \Big(\Delta_\ell \I_\alpha f\Big)(x,t) \Big(\Delta_{\ell-h}\I_\alpha f\Big)^{q-1}(x,t)dxdt
~\leq~\C_{p~q}~ 2^{-\ve h}  \left\| f\right\|_{\L^p(\R^{n+1})}^{q}
 \end{array}
\eeq 
for every $h\in\Z$ nonnegative, provided that $q\in\Z$ is sufficiently large.}
 \v
By applying the lemma to (\ref{Sum I_alpha})-(\ref{Ortho Expansion}), we  have
$\left\|\I_\alpha f\right\|_{\L^{q}(\R^{n+1})}\leq\C_{p~q}\left\|f\right\|_{\L^{p}(\R^{n+1})}$ for $q\in\Z$ sufficiently large.

Observe that $\I_\alpha$ defined in (\ref{I_alpha}) is a self-adjoint operator. Let ${\alpha\over n}={1\over p}-{1\over q}={q-1\over q}-{p-1\over p}$. We have
$\I_\alpha\colon \L^p(\R^{n+1})\mt\L^q(\R^{n+1}) \Longleftrightarrow
\I_\alpha\colon \L^{q\over q-1}(\R^{n+1})\mt\L^{p\over p-1}(\R^{n+1})$.

Let $1<p_i<q_i<\infty,~i=1,2$ satisfy
\bel{formula12}
{1\over p}-{1\over q}~=~{\alpha\over n}~=~{1\over p_1}-{1\over q_1}~=~{1\over p_2}-{1\over q_2}.
\eeq
By choosing $q_1$ and $\left({p_2\over p_2-1}\right)\in\Z$  sufficiently large  and keeping the equalities hold in (\ref{formula12}), we simultaneously have 
 $\left\|\I_\alpha f\right\|_{\L^{q_1}(\R^{n+1})}\leq\C_{p_1~q_1}\left\|f\right\|_{\L^{p_1}(\R^{n+1})},~\left\|\I_\alpha f\right\|_{\L^{q_2}(\R^{n+1})}\leq\C_{p_2,~q_2}\left\|f\right\|_{\L^{p_2}(\R^{n+1})}$ whereas
${1\over p}={1-t\over p_1}+{t\over p_2},~{1\over q}={1-t\over q_1}+{t\over q_2}$ for some $0\leq t\leq 1$. 
By applying Riesz interpolation theorem, we obtain the desired result in (\ref{Result}).

 \section{Pointwise estimates on partial operators}
 \setcounter{equation}{0}
Let $\eta\ge1$ be fixed. For every $0\leq\ell\leq \eta$, we define $\rho_\ell(x,t)\in\R$ implicitly by
\bel{rho}
\left(2^{\rho_\ell(x,t)-\ell} 2^{\rho_\ell(x,t) n}\right)^{1\over p}~=~{\vartheta_\ell^{1\over p}(x,t)\left\| f\right\|_{\L^p(\R^{n+1})}\over \Big(\M_\eta f\Big)(x,t)},\qquad a.e~~(x,t)\in\R^n\times\R.
\eeq
Let $\Lambda_{\ell j}$ and $\Gamma^\nu$ be given in (\ref{Cone}) and (\ref{Gamma}). From (\ref{Inclusion R^v_j})-(\ref{R size}), we have
\bel{Partial Est1 v}
\begin{array}{lr}\ds
\iint_{\Lambda_{\ell j}\cap\Gamma^\nu_\eta}f(x-y,t-s)\left({1\over |s|+|y|}\right)^{n-\alpha}\left({1\over |s|-|y|}\right)^{1-{\alpha\over n}} dyds
\\\\ \ds
~\leq~\C~{1\over 2^{j(n-\alpha)}} {1\over 2^{(j-\ell)(1-{\alpha\over n})}}\iint_{\Lambda_{\ell j}\cap\Gamma^\nu_\eta}f(x-y,t-s)dyds
\\\\ \ds
~=~\C~\left(2^{j-\ell}2^{jn}\right)^{\alpha\over n}{1\over 2^{j-\ell}2^{jn}}\iint_{\Lambda_{\ell j}\cap\Gamma^\nu_\eta}f(x-y,t-s)dyds
\\\\ \ds
~=~\C~\left( 2^{j-\ell}2^{jn}\right)^{\alpha\over n} 2^{-\eta(n-1)} \left\{ {1\over 2^{j-\ell}2^j 2^{(j-\eta)(n-1)}}\iint_{\Lambda_{\ell j}\cap\Gamma^\nu_\eta}f(x-y,t-s)dyds\right\}
\\\\ \ds
~\leq~\C~\left( 2^{j-\ell}2^{jn}\right)^{\alpha\over n} 2^{-\eta(n-1)} \Big(\M^\nu_\eta f\Big)(x,t).
\end{array}
\eeq
Recall the maximal operator $\M_\eta$ defined in (\ref{M*}). We have
\bel{Partial Est1}
\begin{array}{lr}\ds
\iint_{\Lambda_{\ell j}}f(x-y,t-s)\left({1\over |s|+|y|}\right)^{n-\alpha}\left({1\over |s|-|y|}\right)^{1-{\alpha\over n}} dyds
\\\\ \ds
~\leq~\sum_\nu \iint_{\Lambda_{\ell j}\cap\Gamma^\nu_\eta}f(x-y,t-s)\left({1\over |s|+|y|}\right)^{n-\alpha}\left({1\over |s|-|y|}\right)^{1-{\alpha\over n}} dyds
\\\\ \ds
~\leq~\C~\left( 2^{j-\ell}2^{jn}\right)^{\alpha\over n} 2^{-\eta(n-1)} \sum_\nu\Big(\M^\nu_\eta f\Big)(x,t)
\\\\ \ds
~=~\C~\left( 2^{j-\ell}2^{jn}\right)^{\alpha\over n} \Big(\M_\eta f\Big)(x,t).
\end{array}
\eeq
On the other hand, by applying H\"{o}lder inequality, we have
\bel{Partial Est2 Holder}
\begin{array}{lr}\ds
\iint_{\Lambda_{\ell j}}f(x-y,t-s)\left({1\over |s|+|y|}\right)^{n-\alpha}\left({1\over |s|-|y|}\right)^{1-{\alpha\over n}} dyds
\\\\ \ds
~\leq~\left\| f\right\|_{\L^p(\Lambda_{\ell j})}\left\{\iint_{\Lambda_{\ell j}}\left({1\over |s|+|y|}\right)^{(n-\alpha)\left({p\over p-1}\right)}\left({1\over |s|-|y|}\right)^{\left(1-{\alpha\over n}\right)\left({p\over p-1}\right)} dyds\right\}^{p-1\over p}
\\\\ \ds
~\leq~\vartheta_\ell^{1\over p}(x,t)\left\| f\right\|_{\L^p(\R^{n+1})}\left\{\iint_{\Lambda_{\ell j}}\left({1\over |s|+|y|}\right)^{(n-\alpha)\left({p\over p-1}\right)}\left({1\over |s|-|y|}\right)^{\left(1-{\alpha\over n}\right)\left({p\over p-1}\right)} dyds\right\}^{p-1\over p}
\\\\ \ds
~\leq~\C~\vartheta_\ell^{1\over p}(x,t)\left\| f\right\|_{\L^p(\R^{n+1})} \left\{  {1\over 2^{j(n-\alpha)\left({p\over p-1}\right)}} {1\over 2^{(j-\ell)\left(1-{\alpha\over n}\right)\left({p\over p-1}\right)}} 2^{j-\ell}2^{jn}\right\}^{p-1\over p}
\\\\ \ds
~=~\C~\vartheta_\ell^{1\over p}(x,t)\left\| f\right\|_{\L^p(\R^{n+1})} \left\{ 2^{jn\left({p-1\over p}\right)-j(n-\alpha)}2^{(j-\ell)\left({p-1\over p}\right)-(j-\ell)\left(1-{\alpha\over n}\right)}\right\}
\\\\ \ds
~=~\C~\vartheta_\ell^{1\over p}(x,t)\left\| f\right\|_{\L^p(\R^{n+1})} \left\{ 2^{jn\left({\alpha\over n}-1+\left({p-1\over p}\right)\right)}2^{(j-\ell)\left({\alpha\over n}-1+\left({p-1\over p}\right)\right)}\right\}
\\\\ \ds
~=~\C~\vartheta_\ell^{1\over p}(x,t)\left\| f\right\|_{\L^p(\R^{n+1})} \left( 2^{j-\ell}2^{jn}\right)^{{\alpha\over n}-{1\over p}}.
\end{array}
\eeq
By inserting (\ref{rho}) into (\ref{Partial Est1}), for every $\rho_\ell(x,t)<j<\infty$, we have
\bel{Est1 partial}
\begin{array}{lr}\ds
\iint_{\Lambda_{\ell j}}f(x-y,t-s)\left({1\over |s|+|y|}\right)^{n-\alpha}\left({1\over |s|-|y|}\right)^{1-{\alpha\over n}} dyds
\\\\ \ds
~\leq~\C\left( 2^{j-\ell}2^{jn}\right)^{\alpha\over n} \Big(\M_\eta f\Big)(x,t)
\\\\ \ds
~=~\C~2^{\left(j-\rho_\ell(x,t)\right)\left({n+1\over n}\right)\alpha} \left( 2^{\rho_\ell(x,t)-\ell}2^{\rho_\ell(x,t)n}\right)^{\alpha\over n} \Big(\M_\eta f\Big)(x,t)
\\\\ \ds
~=~\C~2^{\left(j-\rho_\ell(x,t)\right)\left({n+1\over n}\right)\alpha} \left( 2^{\rho_\ell(x,t)-\ell}2^{\rho_\ell(x,t)n}\right)^{{1\over p}-{1\over q}} \Big(\M_\eta f\Big)(x,t)
\\\\ \ds
~=~\C~2^{\left(j-\rho_\ell(x,t)\right)\left({n+1\over n}\right)\alpha} \left\{{\vartheta_\ell^{1\over p}(x,t)\left\| f\right\|_{\L^p(\R^{n+1})}\over \Big(\M_\eta f\Big)(x,t)}\right\}^{1-{p\over q}}  \Big(\M_\eta f\Big)(x,t)
\\\\ \ds
~=~\C~2^{\left(j-\rho_\ell(x,t)\right)\left({n+1\over n}\right)\alpha}\vartheta_\ell^{{1\over p}-{1\over q}}(x,t) \Big(\M_\eta f\Big)^{p\over q}(x,t)\left\| f\right\|_{\L^p(\R^{n+1})}^{1-{p\over q}}.
\end{array}
\eeq
By inserting (\ref{rho}) into (\ref{Partial Est2 Holder}), for every $\rho_\ell(x,t)<j<\infty$, we have
\bel{Est2 partial}
\begin{array}{lr}\ds
\iint_{\Lambda_{\ell j}}f(x-y,t-s)\left({1\over |s|+|y|}\right)^{n-\alpha}\left({1\over |s|-|y|}\right)^{1-{\alpha\over n}} dyds
\\\\ \ds
~\leq~\C~\vartheta_\ell^{1\over p}(x,t)\left\| f\right\|_{\L^p(\R^{n+1})} ~\left( 2^{j-\ell}2^{jn}\right)^{{\alpha\over n}-{1\over p}}
\\\\ \ds
~=~\C~2^{-\left(j-\rho_\ell(x,t)\right)(n+1)/q}\vartheta_\ell^{1\over p}(x,t)\left\| f\right\|_{\L^p(\R^{n+1})} \left( 2^{\rho_\ell(x,t)-\ell}2^{\rho_\ell(x,t)n}\right)^{-{1\over q}} \qquad (~{1/ q}>0~)
\\\\ \ds
~=~\C~2^{-\left(j-\rho_\ell(x,t)\right)(n+1)/q}\vartheta_\ell^{1\over p}(x,t)\left\| f\right\|_{\L^p(\R^{n+1})}\left\{{\vartheta_\ell^{1\over p}(x,t)\left\| f\right\|_{\L^p(\R^{n+1})}\over \Big(\M_\eta f\Big)(x,t)}\right\}^{-{p\over q}} 
\\\\ \ds
~=~\C~2^{-\left(j-\rho_\ell(x,t)\right)(n+1)/q}\vartheta_\ell^{{1\over p}-{1\over q}}(x,t)\Big(\M_\eta f\Big)^{p\over q}(x,t) \left\| f\right\|_{\L^p(\R^{n+1})}^{1-{p\over q}}.
\end{array}
\eeq
By putting together (\ref{Est1 partial}) and (\ref{Est2 partial}), we find
\bel{Est regularity partial}
\begin{array}{lr}\ds
\iint_{\Lambda_{\ell j}}f(x-y,t-s)\left({1\over |s|+|y|}\right)^{n-\alpha}\left({1\over |s|-|y|}\right)^{1-{\alpha\over n}} dyds
\\\\ \ds
~\leq~\C~2^{-\left|j-\rho_\ell(x,t)\right|(n+1)\min\left\{{\alpha\over n}, {1\over q}\right\}}\vartheta_\ell^{{1\over p}-{1\over q}}(x,t)\Big(\M_\eta f\Big)^{p\over q}(x,t) \left\| f\right\|_{\L^p(\R^{n+1})}^{1-{p\over q}}.
\end{array}
\eeq
By summing all $j\in\Z$ and using (\ref{Est regularity partial}), we obtain (\ref{Regularity est}) as desired.

\section{Some preliminary results}
\setcounter{equation}{0}
Denote
 \bel{partial lj}
 \Big(\Delta_{\ell j} \I_\alpha f\Big)(x,t) ~=~\iint_{\Lambda_{\ell j}}f(x-y,t-s)\left({1\over |s|+|y|}\right)^{n-\alpha}\left({1\over |s|-|y|}\right)^{1-{\alpha\over n}} dyds 
 \eeq
 for every $\ell\ge0, j\in\Z$.

Let $j,k_i\in\Z, i=1,2,\ldots, q-1$ for $q\in\Z$ sufficiently large. We have
\bel{Ortho Sum}
\begin{array}{lr}\ds
\iint_{\R^{n+1}} \sum_{\ell=h}^\eta \Big(\Delta_\ell \I_\alpha f\Big)(x,t) \Big(\Delta_{\ell-h}\I_\alpha f\Big)^{q-1}(x,t)dxdt
\\\\ \ds
~=~\sum_{j,k_1,k_2,\ldots,k_{q-1}\in\Z}  \iint_{\R^{n+1}}\sum_{\ell=h}^\eta \Big(\Delta_{\ell j} \I_\alpha f\Big)(x,t) \prod_{i=1}^{q-1}\Big(\Delta_{\ell-h~k_i}\I_\alpha f\Big)(x,t)dxdt
\end{array}
\eeq
for every $0\leq h\leq \eta<\infty$.
\vsk

In order to prove (\ref{Ortho Result}), we develop a $3$-fold estimate by considering
\bel{Split}
\sum_{j,k_1,k_2,\ldots,k_{q-1}\in\Z}~=~\sum_{\mathcal{G}_1}+\sum_{\mathcal{G}_2}+\sum_{\mathcal{G}_3}
\eeq
where
\bel{G_1,2,3}
\begin{array}{ccc}
k~=~\min\{k_i,i=1,2,\ldots,q-1\},
\\\\ \ds
\mathcal{G}_1~=~\left\{j,k_1,k_2,\ldots,k_{q-1}\in\Z~\colon~j-h~\ge~k-2\right\},
\\\\ \ds
\mathcal{G}_2~=~\left\{j,k_1,k_2,\ldots,k_{q-1}\in\Z~\colon~j~\leq~k\right\},
\\\\ \ds
\mathcal{G}_3~=~\left\{j,k_1,k_2,\ldots,k_{q-1}\in\Z~\colon~j-h~<~k-2~<~j-2\right\}.
\end{array}
\eeq
By using (\ref{Est regularity partial}), we have
\bel{G est}
\begin{array}{lr}\ds
\Big(\Delta_{\ell j} \I_\alpha f\Big)(x,t) \prod_{i=1}^{q-1} \Big(\Delta_{\ell-h~k_i}\I_\alpha f\Big)(x,t)
\\\\ \ds
~\leq~\C_{p~q}~ 2^{-\left|j-\rho_\ell(x,t)\right|(n+1)\min\left\{{\alpha\over n}, {1\over q}\right\}}  \vartheta^{{1\over p}-{1\over q}}_\ell(x,t)\Big(\M_\eta f\Big)^{p\over q}(x,t)\left\| f\right\|_{\L^p(\R^{n+1})}^{1-{p\over q}}
\\\\ \ds~~~~
\prod_{i=1}^{q-1}2^{-\left|k_i-\rho_{\ell-h}(x,t)\right|(n+1)\min\left\{{\alpha\over n}, {1\over q}\right\}} 
\vartheta^{\left({1\over p}-{1\over q}\right)(q-1)}_{\ell-h}(x,t)\Big(\M_\eta f\Big)^{\left({p\over q}\right)(q-1)}(x,t)\left\| f\right\|_{\L^p(\R^{n+1})}^{\left(1-{p\over q}\right)(q-1)}
\\\\ \ds
~=~\C_{p~q}~ 2^{-\left|j-\rho_\ell(x,t)\right|(n+1)\min\left\{{\alpha\over n}, {1\over q}\right\}} \prod_{i=1}^{q-1}2^{-\left|k_i-\rho_{\ell-h}(x,t)\right|(n+1)\min\left\{{\alpha\over n}, {1\over q}\right\}} 
\\\\ \ds~~~~~~~
 \vartheta^{{1\over p}-{1\over q}}_\ell(x,t)\vartheta^{\left({1\over p}-{1\over q}\right)(q-1)}_{\ell-h}(x,t)\Big(\M_\eta f\Big)^p(x,t)\left\| f\right\|_{\L^p(\R^{n+1})}^{q-p}
\end{array}
\eeq
for every $\ell\leq\eta$.

Let $\jmath(\ell\colon x,t)$ and $\kappa_i(\ell-h\colon x,t), i=1,2,\ldots,q-1$  be defined by
\bel{j kappa}
\begin{array}{ccc}
\jmath(\ell\colon x,t)~=~j-\rho_\ell(x,t),
\\\\ \ds
 \kappa_i(\ell-h\colon x,t)~=~k_i-\rho_{\ell-h}(x,t),~~i=1,2,\ldots,q-1,
\\\\ \ds
\kappa(\ell-h\colon x,t)~=~k-\rho_{\ell-h}(x,t).
\end{array}
\eeq

Recall from (\ref{rho}). We have
\bel{rho rewrite}
\vartheta_\ell(x,t)
~=~{\Big(\M_\eta f\Big)^p(x,t)\over \left\| f\right\|_{\L^p(\R^{n+1})}^p}\left(2^{\rho_\ell(x,t)-\ell} 2^{\rho_\ell(x,t) n}\right).
\eeq

{\bf Case 1}: Let $j-h\ge k-2$. Suppose $\rho_\ell(x,t)-\rho_{\ell-h}(x,t)\ge (1-\sigma)h$ for some $0<\sigma<1$.
We have
\bel{G_1 est1}
\begin{array}{lr}\ds
{\vartheta_{\ell-h}(x,t)\over\vartheta_{\ell}(x,t)}~=~\left( 2^{\rho_\ell(x,t)-\ell} 2^{-\rho_{\ell-h}(x,t)+\ell-h} 2^{\rho_\ell(x,t)n}2^{-\rho_{\ell-h}(x,t)n}\right)^{-1}
\\\\ \ds~~~~~~~~~~~~~~~~
~=~\left( 2^{\left(\rho_\ell(x,t)-\rho_{\ell-h}(x,t)\right)(n+1)} 2^{-h} \right)^{-1}
\\\\ \ds~~~~~~~~~~~~~~~~
~\leq~\left( 2^{(1-\sigma)(n+1)h}2^{-h} \right)^{-1}
\\\\ \ds~~~~~~~~~~~~~~~~
~=~2^{-h\left(n-\sigma(n+1)\right)}.
\end{array}
\eeq
By bringing (\ref{G_1 est1}) back to (\ref{G est}), we have
\bel{G_1 est2}
\begin{array}{lr}\ds
\Big(\Delta_{\ell j} \I_{\alpha} f\Big)(x,t) \prod_{i=1}^{q-1} \Big(\Delta_{\ell-h~k_i}\I_{\alpha} f\Big)(x,t)
\\\\ \ds
~\leq~\C_{p~q}~2^{-\left|j-\rho_\ell(x,t)\right|(n+1)\min\left\{{\alpha\over n}, {1\over q}\right\}} \prod_{i=1}^{q-1}2^{-\left|k_i-\rho_{\ell-h}(x,t)\right|(n+1)\min\left\{{\alpha\over n}, {1\over q}\right\}} 
\\\\ \ds~~~~~~~
 \vartheta^{{1\over p}-{1\over q}}_\ell(x,t)\vartheta^{\left({1\over p}-{1\over q}\right)(q-1)}_{\ell-h}(x,t)\Big(\M_\eta f\Big)^p(x,t)\left\| f\right\|_{\L^p(\R^{n+1})}^{q-p}
\\\\ \ds
~\leq~\C_{p~q}~2^{-\left|\jmath(\ell\colon x,t)\right|(n+1)\min\left\{{\alpha\over n}, {1\over q}\right\}} \prod_{i=1}^{q-1}2^{-\left|\kappa_i(\ell-h\colon x,t)\right|(n+1)\min\left\{{\alpha\over n}, {1\over q}\right\}} 
\\\\ \ds~~~~~~~
2^{-h\left(n-\sigma(n+1)\right)\left({1\over p}-{1\over q}\right)} \vartheta^{2\left({1\over p}-{1\over q}\right)}_{\ell}(x,t)\vartheta^{\left({1\over p}-{1\over q}\right)(q-2)}_{\ell-h}(x,t)\Big(\M_\eta f\Big)^p(x,t)\left\| f\right\|_{\L^p(\R^{n+1})}^{q-p}
\\\\ \ds
~\leq~\C_{p~q}~2^{-\left|\jmath(\ell\colon x,t)\right|(n+1)\min\left\{{\alpha\over n}, {1\over q}\right\}} \prod_{i=1}^{q-1}2^{-\left|\kappa_i(\ell-h\colon x,t)\right|(n+1)\min\left\{{\alpha\over n}, {1\over q}\right\}} 
\\\\ \ds~~~~~~~
2^{-h\left(n-\sigma(n+1)\right)\left({1\over p}-{1\over q}\right)} \vartheta^{\left({1\over p}-{1\over q}\right)(q-2)}_{\ell-h}(x,t)\Big(\M_\eta f\Big)^p(x,t)\left\| f\right\|_{\L^p(\R^{n+1})}^{q-p}.
\end{array}
\eeq
On the other hand, suppose $\rho_\ell(x,t)-\rho_{\ell-h}(x,t)\leq (1-\sigma)h$. Note that 
\bel{G_1 indices}
\begin{array}{lr}\ds
j-h~=~\jmath(\ell\colon x,t)+\rho_\ell(x,t)-h
\\\\ \ds~~~~~~~
~\ge~\kappa(\ell-h\colon x,t)+\rho_{\ell-h}(x,t)-2~=~k-2
\end{array}
\eeq
implies
\bel{G_1 indices est}
\begin{array}{lr}\ds
\jmath(\ell\colon x,t)-\kappa_\imath(\ell-h\colon x,t)~\ge~h-\left( \rho_\ell(x,t)-\rho_{\ell-h}(x,t)\right)-2
\\\\ \ds~~~~~~~~~~~~~~~~~~~~~~~~~~~~~~~~~~~~~
~\ge~h-(1-\sigma)h-2~=~\sigma h-2.
\end{array}
\eeq
By using (\ref{G_1 indices est}), we have
\bel{G_1 est3}
\begin{array}{lr}\ds
 \Big(\Delta_{\ell j} \I_{\alpha} f\Big)(x,t) \prod_{i=1}^{q-1} \Big(\Delta_{\ell-h~k_i}\I_{\alpha} f\Big)(x,t)
 \\\\ \ds
~\leq~\C_{p~q}~2^{-\left|j-\rho_\ell(x,t)\right|(n+1)\min\left\{{\alpha\over n}, {1\over q}\right\}} \prod_{i=1}^{q-1}2^{-\left|k_i-\rho_{\ell-h}(x,t)\right|(n+1)\min\left\{{\alpha\over n}, {1\over q}\right\}} 
\\\\ \ds~~~~~~~
 \vartheta^{{1\over p}-{1\over q}}_\ell(x,t)\vartheta^{\left({1\over p}-{1\over q}\right)(q-1)}_{\ell-h}(x,t)\Big(\M_\eta f\Big)^p(x,t)\left\| f\right\|_{\L^p(\R^{n+1})}^{q-p} 
\\\\ \ds
~\leq~\C_{p~q}~2^{-\left|\jmath(\ell\colon x,t)\right|(n+1)\min\left\{{\alpha\over n}, {1\over q}\right\}} \prod_{i=1}^{q-1}2^{-\left|\kappa_i(\ell-h\colon x,t)\right|(n+1)\min\left\{{\alpha\over n}, {1\over q}\right\}} 
\\\\ \ds~~~~~~~
 \vartheta^{{1\over p}-{1\over q}}_\ell(x,t)\vartheta^{\left({1\over p}-{1\over q}\right)(q-1)}_{\ell-h}(x,t)\Big(\M_\eta f\Big)^p(x,t)\left\| f\right\|_{\L^p(\R^{n+1})}^{q-p}
\\\\ \ds
~\leq~\C_{p~q}~2^{-\left|\jmath(\ell\colon x,t)-\kappa(\ell-h\colon x,t)\right|\left({n+1\over 2}\right)\min\left\{{\alpha\over n}, {1\over q}\right\}} 2^{-\left|\jmath(\ell\colon x,t)\right|\left({n+1\over 2}\right)\min\left\{{\alpha\over n}, {1\over q}\right\}} \prod_{i=1}^{q-1}2^{-\left|\kappa_i(\ell-h\colon x,t)\right|\left({n+1\over 2}\right)\min\left\{\alpha, {1\over q}\right\}} 
\\\\ \ds~~~~~~~
 \vartheta^{{1\over p}-{1\over q}}_\ell(x,t)\vartheta^{\left({1\over p}-{1\over q}\right)(q-1)}_{\ell-h}(x,t)\Big(\M_\eta f\Big)^p(x,t)\left\| f\right\|_{\L^p(\R^{n+1})}^{q-p}
 \\\\ \ds
~\leq~\C_\alpha~2^{-\sigma h\left({n+1\over 2}\right)\min\left\{{\alpha\over n}, {1\over q}\right\}} 2^{-\left|\jmath(\ell\colon x,t)\right|\left({n+1\over 2}\right)\min\left\{{\alpha\over n}, {1\over q}\right\}} \prod_{i=1}^{q-1}2^{-\left|\kappa_i(\ell-h\colon x,t)\right|\left({n+1\over 2}\right)\min\left\{{\alpha\over n}, {1\over q}\right\}} 
\\\\ \ds~~~~~~~
\vartheta^{\left({1\over p}-{1\over q}\right)(q-2)}_{\ell-h}(x,t)\Big(\M_\eta f\Big)^p(x,t)\left\| f\right\|_{\L^p(\R^{n+1})}^{q-p}.
\end{array}
\eeq
By putting together (\ref{G_1 est2}) and (\ref{G_1 est3}) and choosing $\sigma={2n\over 3n+3}$, we find
\bel{G_1 est final}
\begin{array}{lr}\ds
\Big(\Delta_{\ell j} \I_{\alpha} f\Big)(x,t) \prod_{i=1}^{q-1} \Big(\Delta_{\ell-h~k_i}\I_{\alpha} f\Big)(x,t)
\\\\ \ds
~\leq~\C_{p~q}~2^{- h\left({n\over 3}\right)\min\left\{{\alpha\over n}, {1\over q}\right\}}  2^{-\left|\jmath(\ell\colon x,t)\right|\left({n+1\over 2}\right)\min\left\{{\alpha\over n}, {1\over q}\right\}} \prod_{i=1}^{q-1}2^{-\left|\kappa_i(\ell\colon x,t)\right|\left({n+1\over 2}\right)\min\left\{{\alpha\over n}, {1\over q}\right\}} 
\\\\ \ds~~~~~~~
\vartheta^{\left({1\over p}-{1\over q}\right)(q-2)}_{\ell-h}(x,t)\Big(\M_\eta f\Big)^p(x,t)\left\| f\right\|_{\L^p(\R^{n+1})}^{q-p}
\\\\ \ds
~=~\C_{p~q}~2^{- h\left({n\over 3}\right)\min\left\{{\alpha\over n}, {1\over q}\right\}}  2^{-\left|j-\rho_\ell(x,t)\right|\left({n+1\over 2}\right)\min\left\{{\alpha\over n}, {1\over q}\right\}} \prod_{i=1}^{q-1}2^{-\left|k_i-\rho_{\ell-h}(x,t)\right|\left({n+1\over 2}\right)\min\left\{{\alpha\over n}, {1\over q}\right\}} 
\\\\ \ds~~~~~~~
\vartheta^{\left({1\over p}-{1\over q}\right)(q-2)}_{\ell-h}(x,t)\Big(\M_\eta f\Big)^p(x,t)\left\| f\right\|_{\L^p(\R^{n+1})}^{q-p}.
\end{array}
\eeq
By summing over all $0\leq \ell\leq \eta$ and $j, k_1,\ldots, k_{q-1}$ in $\mathcal{G}_1$, we have
 \bel{G_1 est Sum}
\begin{array}{lr}\ds
\sum_{j,k_1,k_2,\ldots,k_{q-1}\in\mathcal{G}_1}  \iint_{\R^{n+1}}\sum_{\ell=h}^\eta\Big(\Delta_{\ell j} \I_{\alpha} f\Big)(x,t) \prod_{i=1}^{q-1} \Big(\Delta_{\ell-h~k_i}\I_{\alpha} f\Big)(x,t)dxdt
\\\\ \ds
~\leq~\C_{p~q}~2^{- h\left({n\over 3}\right)\min\left\{{\alpha\over n}, {1\over q}\right\}}\left\| f\right\|_{\L^p(\R^{n+1})}^{q-p}
\\\\ \ds~~~~~~~~~~~~~~
\iint_{\R^{n+1}} \sum_{j,k_1,\ldots,k_{q-1}\in\mathcal{G}_1} \sum_{\ell=h}^\eta 2^{-\left|j-\rho_\ell(x,t)\right|\left({n+1\over 2}\right)\min\left\{{\alpha\over n}, {1\over q}\right\}} \prod_{i=1}^{q-1}2^{-\left|k_i-\rho_{\ell-h}(x,t)\right|\left({n+1\over 2}\right)\min\left\{{\alpha\over n}, {1\over q}\right\}} 
\\\\ \ds~~~~~~~~~~~~~~
\vartheta^{\left({1\over p}-{1\over q}\right)(q-2)}_{\ell-h}(x,t)\Big(\M_\eta f\Big)^p(x,t) dxdt
\\\\ \ds
~\leq~\C_{p~q}~2^{- h\left({n\over 3}\right)\min\left\{{\alpha\over n}, {1\over q}\right\}}\left\| f\right\|_{\L^p(\R^{n+1})}^{q-p}
\\\\ \ds~~~~~~~~~~~~~~
\iint_{\R^{n+1}}  \sum_{\ell=h}^\eta \left\{\sum_{\jmath,\kappa_1,\ldots,\kappa_{q-1}\in\Z} 2^{-\left|\jmath(\ell\colon x,t)\right|\left({n+1\over 2}\right)\min\left\{{\alpha\over n}, {1\over q}\right\}} \prod_{i=1}^{q-1}2^{-\left|\kappa_i(\ell\colon x,t)\right|\left({n+1\over 2}\right)\min\left\{{\alpha\over n}, {1\over q}\right\}} \right\}
\\\\ \ds~~~~~~~~~~~~~~~
\vartheta^{\left({1\over p}-{1\over q}\right)(q-2)}_{\ell-h}(x,t)\Big(\M_\eta f\Big)^p(x,t) dxdt
\\\\ \ds
~\leq~\C_{p~q}~2^{- h\left({n\over 3}\right)\min\left\{{\alpha\over n}, {1\over q}\right\}}\left\| f\right\|_{\L^p(\R^{n+1})}^{q-p}
\iint_{\R^{n+1}} \left\{ \sum_{\ell=h}^\infty\vartheta^{\left({1\over p}-{1\over q}\right)(q-2)}_{\ell-h}(x,t)\right\}\Big(\M_\eta f\Big)^p(x,t) dxdt
\\\\ \ds
~\leq~\C_{p~q}~2^{- h\left({n\over 3}\right)\min\left\{{\alpha\over n}, {1\over q}\right\}}\left\| f\right\|_{\L^p(\R^{n+1})}^{q-p}\iint_{\R^{n+1}} \Big(\M_\eta f\Big)^p(x,t) dxdt
\\\\ \ds
~\leq~\C_{p~q}~2^{- h\left({n\over 3}\right)\min\left\{{\alpha\over n}, {1\over q}\right\}}\left\| f\right\|_{\L^p(\R^{n+1})}^{q}.
\end{array}
\eeq
\v
{\bf Case 2}: Let $j\leq k$. Suppose $\rho_\ell(x,t)-\rho_{\ell-h}(x,t)\leq \sigma h$ for some $0<\sigma<1$.
We have
\bel{G_2 est1}
\begin{array}{lr}\ds
{\vartheta_{\ell}(x,t)\over\vartheta_{\ell-h}(x,t)}~=~\left( 2^{\rho_\ell(x,t)-\ell} 2^{-\rho_{\ell-h}(x,t)+\ell-h} 2^{\rho_\ell(x,t)n}2^{-\rho_{\ell-h}(x,t)n}\right)
\\\\ \ds~~~~~~~~~~~~~~~~
~=~\left( 2^{\left(\rho_\ell(x,t)-\rho_{\ell-h}(x,t)\right)(n+1)} 2^{-h} \right)
\\\\ \ds~~~~~~~~~~~~~~~~
~\leq~\left( 2^{\sigma(n+1)h}2^{-h} \right)
\\\\ \ds~~~~~~~~~~~~~~~~
~=~2^{-h\left(1-\sigma(n+1)\right)}.
\end{array}
\eeq
By bringing (\ref{G_2 est1}) back to (\ref{G est}), we have
\bel{G_2 est2}
\begin{array}{lr}\ds
\Big(\Delta_{\ell j} \I_{\alpha} f\Big)(x,t) \prod_{i=1}^{q-1} \Big(\Delta_{\ell-h~k_i}\I_{\alpha} f\Big)(x,t)
\\\\ \ds
~\leq~\C_{p~q}~2^{-\left|j-\rho_\ell(x,t)\right|(n+1)\min\left\{{\alpha\over n}, {1\over q}\right\}} \prod_{i=1}^{q-1}2^{-\left|k_i-\rho_{\ell-h}(x,t)\right|(n+1)\min\left\{{\alpha\over n}, {1\over q}\right\}} 
\\\\ \ds~~~~~~~
 \vartheta^{{1\over p}-{1\over q}}_\ell(x,t)\vartheta^{\left({1\over p}-{1\over q}\right)(q-1)}_{\ell-h}(x,t)\Big(\M_\eta f\Big)^p(x,t)\left\| f\right\|_{\L^p(\R^{n+1})}^{q-p}
\\\\ \ds
~\leq~\C_{p~q}~2^{-\left|\jmath(\ell\colon x,t)\right|(n+1)\min\left\{{\alpha\over n}, {1\over q}\right\}} \prod_{i=1}^{q-1}2^{-\left|\kappa_i(\ell-h\colon x,t)\right|(n+1)\min\left\{{\alpha\over n}, {1\over q}\right\}} 
\\\\ \ds~~~~~~~
2^{-h\left(1-\sigma(n+1)\right)\left({1\over p}-{1\over q}\right)} \vartheta^{\left({1\over p}-{1\over q}\right)(q-2)}_{\ell-h}(x,t)\Big(\M_\eta f\Big)^p(x,t)\left\| f\right\|_{\L^p(\R^{n+1})}^{q-p}.
 \end{array}
\eeq
On the other hand, suppose $\rho_\ell(x,t)-\rho_{\ell-h}(x,t)\ge \sigma h$. Note that 
\bel{G_2 indices}
\begin{array}{lr}\ds
j~=~\jmath(\ell\colon x,t)+\rho_\ell(x,t)
\\\\ \ds~
~\leq~\kappa(\ell-h\colon x,t)+\rho_{\ell-h}(x,t)~=~k 
\end{array}
\eeq
implies
\bel{G_2 indices est}
\begin{array}{lr}\ds
\kappa_\imath(\ell-h\colon x,t)-\jmath(\ell\colon x,t)~\ge~\left( \rho_\ell(x,t)-\rho_{\ell-h}(x,t)\right)~\ge~\sigma h.
\end{array}
\eeq
By using (\ref{G_2 indices est}) instead of (\ref{G_1 indices est}),  the same estimate in (\ref{G_1 est3}) can be obtained. 

By putting together  (\ref{G_1 est3}) and (\ref{G_2 est2}) and choosing $\sigma={2\over 3(n+1)}$, we find
 \bel{G_2 est final}
\begin{array}{lr}\ds
\Big(\Delta_{\ell j} \I_{\alpha} f\Big)(x,t) \prod_{i=1}^{q-1} \Big(\Delta_{\ell-h~k_i}\I_{\alpha} f\Big)(x,t)
\\\\ \ds
~\leq~\C_{p~q}~2^{- h\left({1\over 3}\right)\min\left\{{\alpha\over n}, {1\over q}\right\}}  2^{-\left|\jmath(\ell\colon x,t)\right|\left({n+1\over 2}\right)\min\left\{{\alpha\over n}, {1\over q}\right\}} \prod_{i=1}^{q-1}2^{-\left|\kappa_i(\ell\colon x,t)\right|\left({n+1\over 2}\right)\min\left\{{\alpha\over n}, {1\over q}\right\}} 
\\\\ \ds~~~~~~~
\vartheta^{\left({1\over p}-{1\over q}\right)(q-2)}_{\ell-h}(x,t)\Big(\M_\eta f\Big)^p(x,t)\left\| f\right\|_{\L^p(\R^{n+1})}^{q-p}
\\\\ \ds
~=~\C_{p~q}~2^{- h\left({1\over 3}\right)\min\left\{{\alpha\over n}, {1\over q}\right\}}  2^{-\left|j-\rho_\ell(x,t)\right|\left({n+1\over 2}\right)\min\left\{{\alpha\over n}, {1\over q}\right\}} \prod_{i=1}^{q-1}2^{-\left|k_i-\rho_{\ell-h}(x,t)\right|\left({n+1\over 2}\right)\min\left\{{\alpha\over n}, {1\over q}\right\}} 
\\\\ \ds~~~~~~~
\vartheta^{\left({1\over p}-{1\over q}\right)(q-2)}_{\ell-h}(x,t)\Big(\M_\eta f\Big)^p(x,t)\left\| f\right\|_{\L^p(\R^{n+1})}^{q-p}.
\end{array}
\eeq
By summing over all $0\leq \ell\leq \eta$ and $j, k_1,\ldots, k_{q-1}$ in $\mathcal{G}_2$, we have
 \bel{G_2 est Sum}
\begin{array}{lr}\ds
\sum_{j,k_1,\ldots,k_{q-1}\in\mathcal{G}_2}  \iint_{\R^{n+1}}\sum_{\ell=h}^\eta\Big(\Delta_{\ell j} \I_{\alpha} f\Big)(x,t) \prod_{i=1}^{q-1} \Big(\Delta_{\ell-h~k_i}\I_{\alpha} f\Big)(x,t)dxdt
\\\\ \ds
~\leq~\C_{p~q}~2^{- h\left({1\over 3}\right)\min\left\{{\alpha\over n}, {1\over q}\right\}}\left\| f\right\|_{\L^p(\R^{n+1})}^{q-p}
\\\\ \ds~~~~~~~~~~~~~~
\iint_{\R^{n+1}} \sum_{j,k_1,\ldots,k_{q-1}\in\mathcal{G}_2} \sum_{\ell=h}^\eta 2^{-\left|j-\rho_\ell(x,t)\right|\left({n+1\over 2}\right)\min\left\{{\alpha\over n}, {1\over q}\right\}} \prod_{i=1}^{q-1}2^{-\left|k_i-\rho_{\ell-h}(x,t)\right|\left({n+1\over 2}\right)\min\left\{{\alpha\over n}, {1\over q}\right\}} 
\\\\ \ds~~~~~~~~~~~~~~
\vartheta^{\left({1\over p}-{1\over q}\right)(q-2)}_{\ell-h}(x,t)\Big(\M_\eta f\Big)^p(x,t) dxdt
\\\\ \ds
~\leq~\C_{p~q}~2^{- h\left({1\over 3}\right)\min\left\{{\alpha\over n}, {1\over q}\right\}}\left\| f\right\|_{\L^p(\R^{n+1})}^{q-p}
\\\\ \ds~~~~~~~~~~~~~~
\iint_{\R^{n+1}}  \sum_{\ell=h}^\eta \left\{\sum_{\jmath,\kappa_1,\ldots,\kappa_{q-1}\in\Z} 2^{-\left|\jmath(\ell\colon x,t)\right|\left({n+1\over 2}\right)\min\left\{{\alpha\over n}, {1\over q}\right\}} \prod_{i=1}^{q-1}2^{-\left|\kappa_i(\ell\colon x,t)\right|\left({n+1\over 2}\right)\min\left\{{\alpha\over n}, {1\over q}\right\}} \right\}
\\\\ \ds~~~~~~~~~~~~~~~
\vartheta^{\left({1\over p}-{1\over q}\right)(q-2)}_{\ell-h}(x,t)\Big(\M_\eta f\Big)^p(x,t) dxdt
\\\\ \ds
~\leq~\C_{p~q}~2^{- h\left({1\over 3}\right)\min\left\{{\alpha\over n}, {1\over q}\right\}}\left\| f\right\|_{\L^p(\R^{n+1})}^{q-p}
\iint_{\R^{n+1}}  \left\{\sum_{\ell=h}^\infty\vartheta^{\left({1\over p}-{1\over q}\right)(q-2)}_{\ell-h}(x,t)\right\}\Big(\M_\eta f\Big)^p(x,t) dxdt
\\\\ \ds
~\leq~\C_{p~q}~2^{- h\left({1\over 3}\right)\min\left\{{\alpha\over n}, {1\over q}\right\}}\left\| f\right\|_{\L^p(\R^{n+1})}^{q-p}\iint_{\R^{n+1}} \Big(\M_\eta f\Big)^p(x,t) dxdt
\\\\ \ds
~\leq~\C_{p~q}~2^{- h\left({1\over 3}\right)\min\left\{{\alpha\over n}, {1\over q}\right\}}\left\| f\right\|_{\L^p(\R^{n+1})}^{q}.
\end{array}
\eeq

\section{Intersection of dyadic cones}
 \setcounter{equation}{0}
 Recall $\Lambda_{\ell j}$ from (\ref{Cone}). We define
 \bel{Cone (x,t)}
\begin{array}{cc}\ds
\Lambda_\ell(x,t)~=~\Cup_{j\in\mathbb{Z}}~ \Lambda_{\ell j}(x,t),
\\\\ \ds
\Lambda_{\ell j}(x,t)~=~\left\{(y,s)\in\R^n\times\R~\colon 2^j\leq |s-t|+|y-x| < 2^{j+1}, 2^{j-\ell}\leq |s-t|-|y-x|<2^{j-\ell+1}\right\}
\end{array}
\eeq 
and their dyadic variants
 \bel{Cone*}
\begin{array}{cc}\ds
\Lambda_\ell^\ast(x,t)~=~\Cup_{j\in\mathbb{Z}}~ \Lambda_{\ell j}^\ast(x,t),
\\\\ \ds
\Lambda_{\ell j}^\ast(x,t)~=~\left\{(y,s)\in\R^n\times\R~\colon 2^{j-3}\leq |s-t|+|y-x| < 2^{j+3}, 2^{j-\ell-3}\leq |s-t|-|y-x|<2^{j-\ell+3}\right\}.
\end{array}
\eeq 
\begin{remark}
Let $\Lambda_{\ell j}=\Lambda_{\ell j}(0,0)$ and $\Lambda^\ast_{\ell j}=\Lambda^\ast_{\ell j}(0,0)$. Our  estimates given in section 4 and 5 remain to be valid by replacing $\Lambda_{\ell j}$ with $\Lambda^\ast_{\ell j}$.
\end{remark}
Let $j-h<k-2<j-2$. We set $r$ equal to
\bel{r}
0~\leq~r~=~j-k+\ell-h~<~\ell-2.
\eeq
Note that  $j\ge k$, $\ell-h\ge0$ and $j-h<k-2$. 
 \begin{prop} Let $r$ be given in (\ref{r}). We have
\bel{Inter Size Est}
\begin{array}{lr}\ds
\left| \Lambda_{\ell j}(x,t)\cap\Cap_{i=1}^{q-1}\Lambda_{\ell-h~k_i}(x^i,t^i)\right|
~\leq~\C ~2^{j-k-h}~\left| \Lambda_{r j}^\ast(x,t)\cap\Cap_{i=1}^{q-1}\Lambda_{\ell-h~k_i}^\ast(x^i,t^i)\right|
\end{array}
\eeq
for every $(x,t), (x^i,t^i), i=1,2,\ldots,q-1 \in\R^n\times\R$.
\end{prop}

 {\bf Proof}: We essentially consider
 \bel{nonempty}
\Lambda_{\ell j}(x,t)\cap\Cap_{i=1}^{q-1}\Lambda_{\ell-h~k_i}(x^i,t^i)~\neq~\emptyset.
 \eeq 
 Recall $k=\min\{k_i, i=1,2,\ldots,q-1\}$.  For every $(y,s)\in \Lambda_{\ell j}(x,t)\cap\Cap_{i=1}^{q-1}\Lambda_{\ell-h~k_i}(x^i,t^i)$, we denote  $\mathfrak{B}_{k-(l-h)}(y,s)$ to be the ball of radius $2^{k-(l-h)-1}$ centered on $(y,s)$. Observe that by expanding $\Lambda_{\ell-h~k_i}(x^i,t^i)$, $i=1,2,\ldots,q-1$ to their dyadic variants defined in (\ref{Cone*}), we have
 \bel{Ball contained}
 \mathfrak{B}_{k-(l-h)}(y,s)~~~\subset~~~\Cap_{i=1}^{q-1}\Lambda_{\ell-h~k_i}^\ast(x^i,t^i).
 \eeq
Let $(y,s)$ run through the intersection in (\ref{nonempty}). We define
 \bel{U} 
 \mathcal{U}~=~\Cup_{(y,s)\in \Lambda_{\ell j}(x,t)\cap\Cap_{i=1}^{q-1}\Lambda_{\ell-h~k_i}(x^i,t^i)}  ~\mathfrak{B}_{k-l+h}(y,s) 
\eeq
whereas
\bel{U contained}
\mathcal{U}~~~\subset~~~\Cap_{i=1}^{q-1}\Lambda_{\ell-h~k_i}^\ast(x^i,t^i).
\eeq 
 \begin{remark}
$\mathcal{U}$ containing $\Lambda_{\ell j}(x,t)\cap\Cap_{i=1}^{q-1}\Lambda_{\ell-h~k_i}(x^i,t^i)$ gives an expansion of the set by a factor of $2^{k-(\ell-h)-1}$ in all directions.
 \end{remark}
Let $\gamma_{(x,t)}(\bar{y},\bar{s})$ denote the line passing through  a point $(\bar{y},\bar{s})\in\Lambda_{\ell j}(x,t)$ and the $time$-axis where $x-y=0$ such that it is perpendicular to the tangent plane containing the intersection  $\gamma_{(x,t)}(y,s)\cap\{(y,s)\in\R^n\times\R\colon(s-t)=|y-x|\}$. Note that $\gamma_{(x,t)}(y,s)$ is unique for every given $(y,s)\in\Lambda_{\ell j}(x,t)$.

Let $\mathfrak{R}_{\ell j}$ be the smallest set in $\R^{n+1}$ such that
\bel{R_lj}
\begin{array}{cc}\ds
\mathfrak{R}_{\ell j}~~~\supset~~~\Lambda_{\ell j}(x,t)\cap\Cap_{i=1}^{q-1}\Lambda_{\ell-h~k_i}(x^i,t^i)
\\\\ \ds
\hbox{and~~~~ if}~~ (y,s)\in\mathfrak{R}_{\ell j}~~\hbox{then} ~~\gamma_{(x,t)}(y,s)\cap\Lambda_{\ell j}(x,t)~~~\subset~~~\mathfrak{R}_{\ell j}.
\end{array}
\eeq
Let $\mathfrak{R}_{\ell j}^r$ for $j-\ell<j-r-2$ be the smallest set in $\R^{n+1}$ such that
\bel{R_lj^r}
\begin{array}{cc}\ds
\mathfrak{R}_{ \ell  j}^r~~~\supset~~~\Lambda_{\ell j}(x,t)\cap\Cap_{i=1}^{q-1}\Lambda_{\ell-h~k_i}(x^i,t^i)
\\\\ \ds
\hbox{and ~~~~if}~~(y,s)\in\mathfrak{R}_{\ell j}^r~~\hbox{then} ~~\gamma_{(x,t)}(y,s)~\cap~\Cup_{l=\ell}^{r-3} \Lambda_{l j}(x,t)~~~\subset~~~\mathfrak{R}_{\ell j}^r.
\end{array}
\eeq
Recall from (\ref{Ball contained})-(\ref{U contained}) and {\bf Remark 6.2}. By taking into account that 
\bel{distance}
2^{k-(\ell-h)-1}-2^{j-\ell}~=~2^{j-r-1}-2^{j-\ell}~>~2^{j-r-2},
\eeq
we must have 
\bel{Inclusion}
\Lambda_{r-3, j}(x,t)\cap \mathfrak{R}_{\ell j}^r~~~\subset ~~~\Lambda_{r-3, j}(x,t)\cap\mathcal{U}.
\eeq
\begin{figure}[h]
\centering
\includegraphics[scale=0.40]{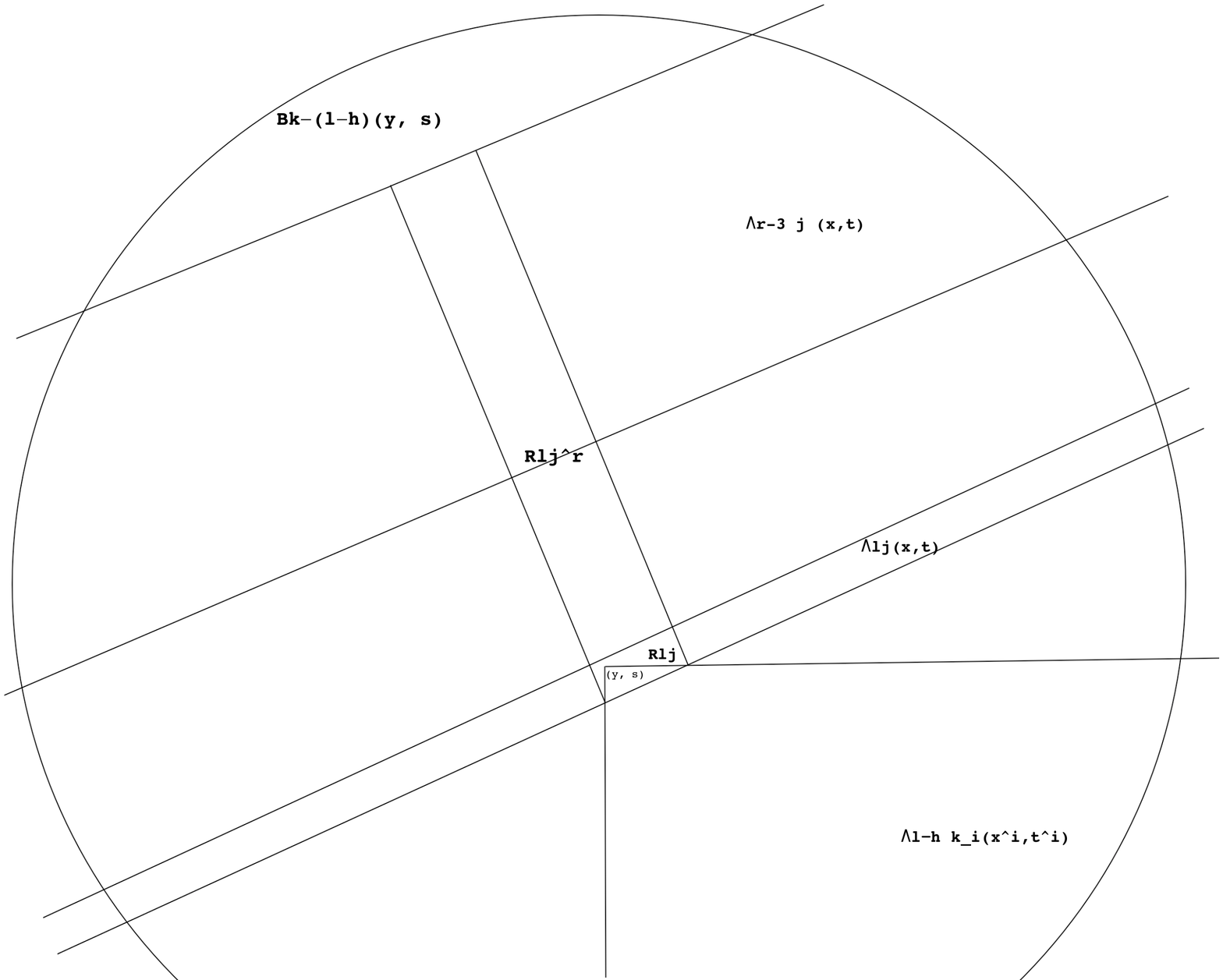}
\end{figure}

By definition of $\mathfrak{R}_{\ell j}$ and $\mathfrak{R}_{\ell j}^r$ defined in (\ref{R_lj}) and (\ref{R_lj^r}), it is a geometric fact that
\bel{computation result}
\begin{array}{lr}\ds
\left|\mathfrak{R}_{\ell j}\right|~\leq~\C~2^{j-r}2^{-(j-\ell)}  \left|\mathfrak{R}_{r\ell j}\right|
\\\\ \ds~~~~~~~
~\leq~\C~2^{\ell-r} \left|\mathfrak{R}_{\ell j}^r\right|
~=~\C~2^{j-k-h} \left|\mathfrak{R}_{\ell j}^r\right|
\\\\ \ds~~~~~~~
~\leq~\C~2^{j-k-h}\left|\Lambda_{r-3, j}(x,t)\cap \mathfrak{R}_{\ell j}^r\right|.
\end{array}
\eeq
By putting together (\ref{Inclusion})-(\ref{computation result}), we find 
 \bel{Lambda_3 dominate}
 \begin{array}{lr}\ds
\left| \Lambda_{\ell j}(x,t)\cap\Cap_{i=1}^{q-1}\Lambda_{\ell-h~k_i}(x^i,t^i)\right|~\leq~\left|\mathfrak{R}_{\ell j}\right|~\leq~\C ~2^{j-k-h}\left| \Lambda_{r-3~ j}(x,t)\cap\mathcal{U} \right|.
\end{array}
 \eeq
 Lastly, by definition of $\Lambda^\ast_{\ell j}(x,t)$ in (\ref{Cone*}), we have
 \bel{Lambda_3 inclusion}
 \Lambda_{r-3~j}(x,t)~\subset~\Lambda^\ast_{rj}(x,t).
 \eeq
 Together with the inclusion given in (\ref{U}), we obtain (\ref{Inter Size Est}) as required.
 \endproof

\section{Proof of almost orthogonality}
\setcounter{equation}{0}
For every $(x,t)\in\Lambda$ in (\ref{I_alpha}), we denote 
\bel{Omega^alpha}
\Omega^\alpha(x,t)~=~\left({1\over |t|+|x|}\right)^{n-\alpha}\left({1\over |t|-|x|}\right)^{1-{\alpha\over n}}.
\eeq
Let $j-h\leq k-2<j-2$. From direct computation, we have
\bel{I_alpha^q v}
\begin{array}{lr}\ds
\iint_{\R^{n+1}} \Big(\Delta_{\ell j}\I_\alpha f\Big)(x,t)\prod_{i=1}^{q-1}\Big(\Delta_{\ell-h~k_i}\I_\alpha f\Big)(x,t)dxdt
\\\\ \ds
~=~ \iint\cdots\iint_{\R^{n+1}\times\cdots\times\R^{n+1}} f(y,s)\prod_{i=1}^{q-1}f(y^i,s^i) 
\\\\ \ds~~~~~~~~
\left\{\iint_{\Lambda_{\ell j}(y,s)\cap\Cap_{i=1}^{q-1}\Lambda_{\ell-h~k_i}(y^i,s^i)}\Omega^\alpha(x-y, t-s)\prod_{i=1}^{q-1}\Omega^\alpha(x-y^i,t-s^i)dxdt  \right\}dyds\prod_{i=1}^{q-1}dy^ids^i
\end{array}
\eeq
where $\Lambda_{\ell j}(y,s)$ is defined in (\ref{Cone (x,t)}).

Set $r=j-k+\ell-h$ as in (\ref{r}). Recall $\Lambda^\ast_{\ell j}(x,t)$ defined in (\ref{Cone*}). By applying {\bf Proposition 6.1}, we have
\bel{Crucial Est}
\begin{array}{lr}\ds
\iint_{\Lambda_{\ell j}(y,s)\cap\Cap_{i=1}^{q-1}\Lambda_{\ell-h~k_i}(y^i,s^i)}\Omega^\alpha(x-y, t-s)\prod_{i=1}^{q-1}\Omega^\alpha(x-y^i,t-s^i)dxdt 
\\\\ \ds
~\leq~\C~\left({1\over 2^j}\right)^{n-\alpha}\left({1\over 2^{j-\ell}}\right)^{1-{\alpha\over n}}\prod_{i=1}^{q-1}
\left({1\over 2^{k_i}}\right)^{n-\alpha}\left({1\over 2^{k_i-(\ell-h)}}\right)^{1-{\alpha\over n}}
\left|\Lambda_{\ell j}(y,s)\cap\Cap_{i=1}^{q-1}\Lambda_{\ell-h~k_i}(y^i,s^i)\right|
\\\\ \ds
~=~\C~2^{(j-k-h)\left({\alpha\over n}-1\right)}\left({1\over 2^j}\right)^{n-\alpha}\left({1\over 2^{j-r}}\right)^{1-{\alpha\over n}}\prod_{i=1}^{q-1}
\left({1\over 2^{k_i}}\right)^{n-\alpha}\left({1\over 2^{k_i-(\ell-h)}}\right)^{1-{\alpha\over n}}
\left|\Lambda_{\ell j}(y,s)\cap\Cap_{i=1}^{q-1}\Lambda_{\ell-h~k_i}(y^i,s^i)\right|
\\\\ \ds
~\leq~\C~2^{(j-k-h)\left({\alpha\over n}-1\right)}\left({1\over 2^j}\right)^{n-\alpha}\left({1\over 2^{j-r}}\right)^{1-{\alpha\over n}}\prod_{i=1}^{q-1}
\left({1\over 2^{k_i}}\right)^{n-\alpha}\left({1\over 2^{k_i-(\ell-h)}}\right)^{1-{\alpha\over n}}
\\\\ \ds~~~~~~~
2^{j-k-h}~\left| \Lambda_{r j}^\ast(y,s)\cap\Cap_{i=1}^{q-1}\Lambda_{\ell-h~k_i}^\ast(y^i,s^i)\right|\qquad \hbox{\small{by (\ref{Inter Size Est})}}
\\\\ \ds
~=~\C~2^{(j-k-h)\left({\alpha\over n}\right)}\left({1\over 2^j}\right)^{n-\alpha}\left({1\over 2^{j-r}}\right)^{1-{\alpha\over n}}\prod_{i=1}^{q-1}
\left({1\over 2^{k_i}}\right)^{n-\alpha}\left({1\over 2^{k_i-(\ell-h)}}\right)^{1-{\alpha\over n}}
\left| \Lambda_{r j}^\ast(y,s)\cap\Cap_{i=1}^{q-1}\Lambda_{\ell-h~k_i}^\ast(y^i,s^i)\right|
\\\\ \ds
~\leq~\C~2^{(j-k-h)\left({\alpha\over n}\right)}\iint_{\Lambda_{r j}^\ast(y,s)\cap\Gamma^\nu_\eta\cap\Cap_{i=1}^{q-1}\Lambda_{\ell-h~k_i}^\ast(y^i,s^i)}\Omega^\alpha(x-y, t-s)\prod_{i=1}^{q-1}\Omega^\alpha(x-y^i,t-s^i)dxdt. 
\end{array}
\eeq
By bringing (\ref{Crucial Est}) back to (\ref{I_alpha^q v}), we have
\bel{I_alpha^q v Est Sum}
\begin{array}{lr}\ds
\iint_{\R^{n+1}} \Big(\Delta_{\ell j}\I_\alpha f\Big)(x,t)\prod_{i=1}^{q-1}\Big(\Delta_{\ell-h~k_i}\I_\alpha f\Big)(x,t)dxdt
\\\\ \ds
~=~  \iint\cdots\iint_{\R^{n+1}\times\cdots\times\R^{n+1}} f(y,s)\prod_{i=1}^{q-1}f(y^i,s^i) 
\\\\ \ds~~~~~~~~
\left\{\iint_{\Lambda_{\ell j}(y,s)\cap\Cap_{i=1}^{q-1}\Lambda_{\ell-h~k_i}(y^i,s^i)}\Omega^\alpha(x-y, t-s)\prod_{i=1}^{q-1}\Omega^\alpha(x-y^i,t-s^i)dxdt  \right\}dyds\prod_{i=1}^{q-1}dy^ids^i
\\\\ \ds
~\leq~ \C~2^{(j-k-h)\left({\alpha\over n}\right)} \iint\cdots\iint_{\R^{n+1}\times\cdots\times\R^{n+1}} f(y,s)\prod_{i=1}^{q-1}f(y^i,s^i) 
\\\\ \ds~~~~~~~~
\left\{\iint_{\Lambda_{r j}^\ast(x,t)\cap\Cap_{i=1}^{q-1}\Lambda_{\ell-h~k_i}^\ast(x^i,t^i)}\Omega^\alpha(x-y, t-s)\prod_{i=1}^{q-1}\Omega^\alpha(x-y^i,t-s^i)dxdt  \right\}dyds\prod_{i=1}^{q-1}dy^ids^i
\\\\ \ds
~=~ \C~2^{(j-k-h)\left({\alpha\over n}\right)}~\iint_{\R^{n+1}} \left\{\iint_{\Lambda_{r j}^\ast}f(x-y,t-s)\Omega^\alpha(y,s) dyds\right\}
\\\\ \ds~~~~~~~~~~~~~~~~~~~~~~~
\prod_{i=1}^{q-1} \left\{\iint_{\Lambda_{\ell-h~k_i}^\ast}f(x-y^i,t-s^i)\Omega^\alpha(y^i,s^i)dy^ids^i\right\}dxdt.
\end{array}
\eeq
Observe that $j-r=k-(\ell-h)$ ( or $j-(r-(\ell-h))=k$ ) where $r-\ell+h=j-k>0$. We are now back to the situition in {\bf Case 1} for which $\ell$ and $h$ are replaced by $r$ and $j-k$ respectively. 

Recall from {\bf Remark 6.1}. The estimate in  (\ref{G_1 est final}) remains valid with $\Lambda_{\ell j}, \Lambda_{\ell-h~k_i}, i=1,2,\ldots,q-1$ replaced by $\Lambda^*_{\ell j}, \Lambda^*_{\ell-h~k_i}, i=1,2,\ldots,q-1$.
By applying the resulting estimate to (\ref{I_alpha^q v Est Sum}), we have
\bel{G_3 est}
\begin{array}{lr}\ds
\iint_{\R^{n+1}} \Big(\Delta_{\ell j}\I_\alpha f\Big)(x,t)\prod_{i=1}^{q-1}\Big(\Delta_{\ell-h~k_i}\I_\alpha f\Big)(x,t)dxdt
\\\\ \ds
~\leq~ \C~2^{(j-k-h)\left({\alpha\over n}\right)}~\iint_{\R^{n+1}} \left\{\iint_{\Lambda_{r j}^\ast}f(x-y,t-s)\Omega^\alpha(y,s) dyds\right\}
\\\\ \ds~~~~~~~~~~~~~~~~~~~~~~~~~~
\prod_{i=1}^{q-1} \left\{\iint_{\Lambda_{\ell-h~k_i}^\ast}f(x-y^i,t-s^i)\Omega^\alpha(y^i,s^i)dy^ids^i\right\}dxdt
\\\\ \ds
~\leq~\C~2^{(j-k-h)\left({\alpha\over n}\right)}

2^{- (j-k)\left({n\over 3}\right)\min\left\{{\alpha\over n}, {1\over q}\right\}}  2^{-\left|j-\rho_r(x,t)\right|\left({n+1\over 2}\right)\min\left\{{\alpha\over n}, {1\over q}\right\}} \prod_{i=1}^{q-1}2^{-\left|k_i-\rho_{\ell-h}(x,t)\right|\left({n+1\over 2}\right)\min\left\{{\alpha\over n}, {1\over q}\right\}} 
\\\\ \ds~~~~~~~
\vartheta^{\left({1\over p}-{1\over q}\right)(q-2)}_{\ell-h}(x,t)\Big(\M_\eta f\Big)^p(x,t)\left\| f\right\|_{\L^p(\R^{n+1})}^{q-p}.
\end{array}
\eeq
Note that $j-k-h<0$ and $j-k>0$. The last line of (\ref{G_3 est}) can be further bounded by
\bel{G_3 est final}
\begin{array}{lr}\ds
\C~2^{(j-k-h)\left({1\over 3}\right)\min\left\{{\alpha\over n},{1\over q}\right\}}

2^{- (j-k)\left({2\over 3}\right)\min\left\{{\alpha\over n}, {1\over q}\right\}}  2^{-\left|j-\rho_r(x,t)\right|\left({n+1\over 2}\right)\min\left\{{\alpha\over n}, {1\over q}\right\}} \prod_{i=1}^{q-1}2^{-\left|k_i-\rho_{\ell-h}(x,t)\right|\left({n+1\over 2}\right)\min\left\{{\alpha\over n}, {1\over q}\right\}} 
\\\\ \ds~~~~~~~
\vartheta^{\left({1\over p}-{1\over q}\right)(q-2)}_{\ell-h}(x,t)\Big(\M_\eta f\Big)^p(x,t)\left\| f\right\|_{\L^p(\R^{n+1})}^{q-p}\qquad ( n\ge2 )
\\\\ \ds
~\leq~\C~2^{-h\left({1\over 3}\right)\min\left\{{\alpha\over n},{1\over q}\right\}}

2^{- (j-k)\left({1\over 3}\right)\min\left\{{\alpha\over n}, {1\over q}\right\}}   \prod_{i=1}^{q-1}2^{-\left|k_i-\rho_{\ell-h}(x,t)\right|\left({n+1\over 2}\right)\min\left\{{\alpha\over n}, {1\over q}\right\}} 
\\\\ \ds~~~~~~~
\vartheta^{\left({1\over p}-{1\over q}\right)(q-2)}_{\ell-h}(x,t)\Big(\M_\eta f\Big)^p(x,t)\left\| f\right\|_{\L^p(\R^{n+1})}^{q-p}.
\end{array}
\eeq
Let $m=j-k$. By summing over all $0\leq \ell\leq \eta$ and $j, k_1,\ldots, k_{q-1}\in\mathcal{G}_3$ given in (\ref{G_1,2,3}), we have
 \bel{G_3 est Sum}
\begin{array}{lr}\ds
\sum_{j,k_1,\ldots,k_{q-1}\in\mathcal{G}_3}  \iint_{\R^{n+1}}\sum_{\ell=h}^\eta\Big(\Delta_{\ell j} \I_{\alpha} f\Big)(x,t) \prod_{i=1}^{q-1} \Big(\Delta_{\ell-h~k_i}\I_{\alpha} f\Big)(x,t)dxdt
\\\\ \ds
~\leq~\C_{p~q}~2^{- h\left({1\over 3}\right)\min\left\{{\alpha\over n}, {1\over q}\right\}}\left\| f\right\|_{\L^p(\R^{n+1})}^{q-p} \sum_{j,k_1,\ldots,k_{q-1}\in\mathcal{G}_3} \sum_{\ell=h}^\eta
\\\\ \ds
\iint_{\R^{n+1}}  2^{-(j-k)\left({1\over 3}\right)\min\left\{{\alpha\over n}, {1\over q}\right\}} \prod_{i=1}^{q-1}2^{-\left|k_i-\rho_{\ell-h}(x,t)\right|\left({n+1\over 2}\right)\min\left\{{\alpha\over n}, {1\over q}\right\}} 
\vartheta^{\left({1\over p}-{1\over q}\right)(q-2)}_{\ell-h}(x,t)\Big(\M_\eta f\Big)^p(x,t) dxdt
\\\\ \ds
~\leq~\C_{p~q}~2^{- h\left({1\over 3}\right)\min\left\{{\alpha\over n}, {1\over q}\right\}}\left\| f\right\|_{\L^p(\R^{n+1})}^{q-p}
\\\\ \ds~~~~~~~~~~~~~~
\iint_{\R^{n+1}}  \sum_{\ell=h}^\eta \left\{~\sum_{m=0}^\infty~\sum_{\kappa_1,\ldots,\kappa_{q-1}\in\Z} 2^{-m\left({1\over 3}\right)\min\left\{{\alpha\over n}, {1\over q}\right\}}\prod_{i=1}^{q-1}2^{-\left|\kappa_i(\ell\colon x,t)\right|\left({n+1\over 2}\right)\min\left\{{\alpha\over n}, {1\over q}\right\}}~ \right\}
\\\\ \ds~~~~~~~~~~~~~~~
\vartheta^{\left({1\over p}-{1\over q}\right)(q-2)}_{\ell-h}(x,t)\Big(\M_\eta f\Big)^p(x,t) dxdt
\\\\ \ds
~\leq~\C_{p~q}~2^{- h\left({1\over 3}\right)\min\left\{{\alpha\over n}, {1\over q}\right\}}\left\| f\right\|_{\L^p(\R^{n+1})}^{q-p}
\iint_{\R^{n+1}}  \left\{\sum_{\ell=h}^\infty\vartheta^{\left({1\over p}-{1\over q}\right)(q-2)}_{\ell-h}(x,t)\right\}\Big(\M_\eta f\Big)^p(x,t) dxdt
\\\\ \ds
~\leq~\C_{p~q}~2^{- h\left({1\over 3}\right)\min\left\{{\alpha\over n}, {1\over q}\right\}}\left\| f\right\|_{\L^p(\R^{n+1})}^{q-p}\iint_{\R^{n+1}} \Big(\M_\eta f\Big)^p(x,t) dxdt
\\\\ \ds
~\leq~\C_{p~q}~2^{- h\left({1\over 3}\right)\min\left\{{\alpha\over n}, {1\over q}\right\}}\left\| f\right\|_{\L^p(\R^{n+1})}^{q}.
\end{array}
\eeq
Lastly, our estimates in (\ref{G_1 est Sum}), (\ref{G_2 est Sum}) and (\ref{G_3 est Sum}) hold for every $1\leq\eta<\infty$. By letting $\eta\mt\infty$, we obtain the desired estimate in (\ref{Ortho Result}).

wangzipeng@westlake.edu.cn

\end{document}